# Finite / Countable State Space Stochastic Processes : Point Processes: Characterization of Associated Auto-Correlation Functions:

Garimella Rama Murthy,
Associate Professor,
IIIT—Hyderabad, Gachibowli,
HYDERABAD-32, AP, INDIA

**ABSTRACT**

In this research paper, the relationship between finite / countable state space stochastic processes and point processes is explored. Utilizing the known relationship between Poisson processes and continuous time Markov chains, finite / countable state space random processes are related to continuous time Markov Chains. Based on the known results for binary random processes, characterization of auto-correlation function of finite state space random processes is explored. An important characterization of corner positive definite matrices is provided.

**1. Introduction:**

Most natural as well as artificial dynamic phenomena are endowed with a stochastic rather than deterministic description. Thus, such a phenomena are mathematically modeled using a stochastic process. Among the class of stochastic processes utilized in stochastic modeling, Markovian processes assume special importance. The main reason is the existence of equilibrium probabilistic behaviour which can be computed efficiently. Thus, many equilibrium performance measures are efficiently determined. The phenomena modeled by Markovian stochastic processes satisfy the property that conditional on present and past, the future depends only on the present. This modeling assumption is satisfied in many stochastic models of natural and artificial phenomena. Phenomena occuring in higher dimensions are modeled using Markov Random Fields.

In some applications, events occur on the real line at random times. Such phenomena are modeled using stochastic processes called point processes. Point processes are utilized in stochastic models arising in queueing, inventory and communication theories. Among point processes, Poisson process is widely utilized in tractable stochastic models. In the case of Poisson process, the interarrival time distributions are exponential. The exponential distribution is the only such distribution which satisfies a property called "memorylessness" [Ross].

Researchers attempted and succeeded in relating Poisson process and continuous time Markov chains. Based on this connection, the author became interested in exploring the relationship between arbitrary finite / countable state space stochastic processes and Point processes. The author hopes that detailed investigation of this relationship enables the analysis of "structured finite / countable state space stochastic processes" using the theory of continuous time Markov chains.

This research paper is organized as follows. In Section 2, relationship between finite / countable state space stochastic processes and point processes is explored. In Section 3, countable state space random processes and countable state space Markov chains are related. In Section 4, characterization of autocorrelation function of finite state space random processes is explored. The research paper concludes in Section 5.

## 2. Relationship between Finite / Countable State Space Stochastic Processes and Point Processes:

Let the state space of a finite state space stochastic process be $\{1, 2, \ldots, N\}$. Considering the random variables representing the successive sojourn times in state "j" ($1 \leq j \leq N$) be

$$T_1^{(j)}, T_2^{(j)}, \ldots\ldots\ldots$$

These random variables can be treated as the inter-arrival times of the "j"th point process. Thus, we have the following claim.

**Claim 1:** An arbitrary finite state space random process can be treated as the superposition of finitely many competing point processes.

The above claim is a generalization of the following claim.

**Claim 2:** Any finite state space Continuous Time Markov Chain (CTMC) can be treated as the superposition of finitely many competing Poisson processes [Cin].

The above claims have natural generalization to countable state space stochastic processes (under some regularity conditions).

Application of the above Claims:

In a large packet network, at routers / switches, one observes the superposition of packet streams. Thus, we are naturally led to the analytical investigation of superposition of point processes (modeling the packet streams). Also, there are several other applications where it is necessary to characterize the superposition of finitely / countably many point processes.

## 3. Countable State Space Random Processes : Countable State Space Markov Chains:

From theoretical as well as practical considerations, one is interested in superposition of countably many point processes. The following claim is relevant in such a discussion and follows from a generalization of the ideas in Section 2.

**Claim 3:** A countable state space stochastic process can be treated as the superposition of countably many point processes.

Now, we are interested when a countable state space stochastic process constitutes a Continuous Time Markov Chain (CTMC) or more generally when it can be approximated by a CTMC. The following claim provides one possible answer in the proper direction and it is well understood in the point process literature. For mathematically formal statement, one can refer to any standard textbook on point processes.

**Claim 4:** The superposition of countably many, "uniformly sparse" point processes constitutes a Poisson process.

The condition of "uniform sparseness" of point processes is well understood, formalized and documented in the literature.

The claim 4 is carefully interpreted to arrive at the relationship between countable state space stochastic processes and countable state space Markov chains.

- Consider an infinite state space stochastic process as the superposition of countably many competing point processes. When these competing point processes are "uniformly sparse", let the stochastic process be of "SPECIAL TYPE".

**Theorem 1:** Any "special type" stochastic process constitutes a countable state space "Uniformized" Continuous Time Markov Chain.

**Remark:** It is well known that the Central Limit Theorem has theoretical as well as practical significance. We expect the above theorem to be of utility in various applications including characterization of packet streams in point-to-point networks.

4. **Characterization of Auto-Correlation Function of Finite State Space Random Processes:**

It is well known from Bochner-Herglotz's Theorem that the autocorrelation function of a stationary random process needs to be a positive definite function. But researchers such as Mc-Millan became interested in characterizing the autocorrelation function ( equals the auto-covariance function when the mean of the stationary random process is zero ) of stationary random process which assumes the values { +1, -1 } ( with equal probability so that the mean is zero ).

$$\rho(n) = E\{X_m X_{m+n}\}$$

Let the class of functions be denoted by $U$. In fact Mc-Millan's stated the following Theorem:

**Theorem 2: ( B. Mc Millan ):**

$$\rho \in U \Leftrightarrow (\text{if and only if})$$
$$\rho(0) = 1 \quad and$$
$$\sum_{n=1}^{N}\sum_{m=1}^{N} \rho(m-n) A_{mn} \geq 0$$

for every N, and all corner-positive matrices $\{A_{mn}\}$, $m, n = 1, 2, \ldots, N$, where a matrix $\{A_{mn}\}$ is corner-positive if

$$\sum_{n=1}^{N}\sum_{m=1}^{N} A_{mn} \epsilon_m \epsilon_n \geq 0$$

for every sequence $(\epsilon_1, \ldots, \epsilon_N)$ with $\epsilon_i = \mp 1$, $i = 1, 2, \ldots, N$.

The above characterization of auto-correlation function of a unit random process can be reformulated as follows:

$$\rho(0) = 1 \quad and$$

**Trace ( R X ) = Trace ( X R ) ≥ 0 ( linear form in the components of X )**

where R is N x N autocorrelation matrix ( corresponding to various lags ) for every N, and all corner-positive matrices $\{X_{mn}\}$, $m, n = 1, 2, \ldots, N$, where a matrix $\{X_{mn}\}$ is corner-positive if

$$\epsilon^T X \epsilon \geq 0 \quad (i.e.\ quadratic\ form\ associated\ with\ X)$$

for every vector $\epsilon^T = (\epsilon_1, \ldots, \epsilon_N)$ lying on the symmetric binary hypercube.

- Thus, we can denote the class U of unit autocorrelation functions as "corner-positive definite" functions. The above characterization problem naturally deals with determining the properties of "corner positive ( definite )" matrices.

As a natural generalization, we introduce the following definitions related to positive definiteness of matrices on interesting subsets of $R^n/C^n$.

(A) V-positive Definite Matrix: If a matrix is positive definite at all points of V i.e.
$X^T A X > 0 \quad for\ all\ X \in V$,
then, the matrix A is said to be a V positive Definite matrix.

Note: The set V can be any finite /countable / uncountable set of $R^n$ or $C^n$.
For instance, if V happens to be bounded, symmetric lattice, we call A to be "Bounded Lattice Positive Definite" matrix. Also V can be the entire lattice in $R^n$.

Note: As a natural generalization, we can define V-positive Definite functions ( in the spirit of Bochner—Herglotz's Theorem ). These could be called Corner Positive Definite functions, Lattice positive definite functions etc.

Now, we start with a characterization of Corner Positive Definite matrices. It should be kept in mind that the basic approach works for bounded-lattice positive definite matrices.

**Theorem 3:** A matrix C ( defined on the field of real or complex numbers ) is Corner Positive Definite (CPD)
( i.e. $X^T C X \geq 0\ for\ all\ vectors\ X\ lying\ on\ the\ symmetric\ binary\ hypercube$ )
if and only if
$X^T E X \geq Trace\ (C)\ for\ all\ vectors\ X\ such\ that$
$X = -sign\ (E X)\ where\ sign\ (.)\ is\ the\ signum\ function$
and E is the symmetric matrix with zero diagonal elements obtained from D, where $D = \frac{1}{2}\ (C + C^T)$.

**Proof:** The proof follows based on non-linear optimization techniques ( maximum principle or Kuhn-Tucker conditions ) and is avoided for brevity. It is related to the Theorem proved in [Rama]. Q.E.D.

**Remark**: Hopfield neural network leads to the optimization of a quadratic form on the symmetric binary hypercube, H. It is proved that the local/global optimum occur at the stable states i.e.
$X \in H\ such\ that\quad X = Sign\ (E X)$.

The author alongwith his student Praveen generalized the Theorem to the case where C is a complex valued matrix [Rama].

- It should be noted that in the above Theorem, we are interested in determining the local minimum of quadratic form evaluated on the symmetric binary hypercube. Let the local minimum vectors on the symmetric binary hypercube be called the "anti-stable states". The ith component of a candidate optimal solution (i.e. an anti-stable state) can be computed by the following successive approximation procedure:
$$X_i(n+1) = -\text{Sign}\left(\sum_{j=1}^{m} E_{ij} X_j(n)\right).$$
As in the case of serial mode of operation of Hopfield neural network, the proof of convergence of the above successive approximation scheme can easily be provided. Details are avoided for brevity. Also the author formulated and solved the problem of minimizing the quadratic form associated with a Hermitian matrix on the hypercube.

- In view of the above discussion, there is no loss of generality in assuming X to be a symmetric zero diagonal matrix in Theorem 3.
- From basic linear algebra, we infer the following:
  Trace (RX) = Trace (XR) with Trace (X) = 0 and X, R are both symmetric matrices. Also, it is well known that both the matrices XR, RX have the same characterisitic polynomial.

Now, we consider the generalization of Mc-Millan's Theorem to the characterization of autocorrelation function (auto-covariance function of a stationary zero mean process) of stationary finite state space random processes.

- Let the class of such autocorrelation functions of a finite state space random process, $X_n$
$$(X_n = \pm j, \quad j \in \{\pm 1, \pm 2, \ldots, \pm M\}$$
with the condition that
$$P\{X_n = +j\} = P\{X_n = -j\} = \frac{1}{2M}$$
be denoted by V.

**Theorem 4:** $\rho \in V \Leftrightarrow$ (if and only if)
$$\rho(0) = \frac{(M+1)(2M+1)}{6} \quad \text{and}$$
$$\sum_{n=1}^{M}\sum_{m=1}^{M} \rho(m-n) B_{mn} \geq 0$$

for every M, and all bounded, symmetric lattice positive definite matrices $\{B_{mn}\}$, $m, n = 1, 2, \ldots, M$, where a matrix $\{B_{mn}\}$ is bounded symmetric lattice positive definite if
$$\sum_{n=1}^{M}\sum_{m=1}^{M} B_{mn} \epsilon_m \epsilon_n \geq 0$$
for every sequence $(\epsilon_1, \ldots, \epsilon_M)$ lying on the bounded symmetric lattice.

**Proof**: Necessity and sufficiency follow from a generalization of the argument used to prove the Mc Millan's Theorem. It is avoided for brevity Q.E.D.

**Note:** The generalization of above Theorem for countable state space random processes follows by doing the optimization on the lattice. It can be compared with the Bochner-Herglotz's Theorem.

**Remark:** Bruck et.al considered the problem of maximization of higher degree forms (than quadratic forms) on the bounded lattice [Rama]. These results naturally suggest the minimization of quadratic forms over the bounded symmetric lattice. Thus, the characterization of bounded symmetric lattice positive definite matrices follows from the results in constrained optimization.

## 5. Conclusions:

In this research paper, point processes are related to finite / countable state space stochastic processes. Using the well known relationship between, Poisson processes and finite state space continuous time Markov chains, structured countable state space random processes are related to countable state space continuous time Markov chain. Also, continuous time Markov chain model of channels in point-to-point networks is proposed.